\input amstex
\input amsppt.sty
\magnification=\magstep1 \hsize=32truecc
 \vsize=22.5truecm
\baselineskip=14truept \NoBlackBoxes
\def\q{\quad}
\def\qq{\qquad}

\def\t{\text}
\def\qtq#1{\q\t{#1}\q}
\def\f{\frac}

\def\a{\alpha}
\def\b{\beta}
\def\sn{k_1+2k_2+\cdots+nk_n}
\def\kn1{k_1!\cdots k_n!}

\def\bi{\binom}

\def\({\left(}
\def\){\right)}
\par\q  Preprint: March 19, 2009
\par\q

\def\Ls#1#2{\Big(\f{#1}{#2}\Big)}

\let \pro=\proclaim
\let \endpro=\endproclaim

\topmatter
\title Some inversion formulas and formulas for Stirling numbers
\endtitle
\author ZHI-Hong Sun\endauthor
\affil School of the Mathematical Sciences, Huaiyin Normal
University,
\\ Huaian, Jiangsu 223001, PR China
\\ E-mail: zhsun$\@$hytc.edu.cn
\\ Homepage: http://www.hytc.edu.cn/xsjl/szh
\endaffil

 \nologo \NoRunningHeads
 \abstract{In
the paper we present some new inversion formulas and two new
formulas for Stirling numbers.
\par\q
}
 \endabstract
 \footnote"" {The author is
supported by the Natural Sciences Foundation of China (grant No.
10971078).}

\endtopmatter
\document

\subheading{1. Introduction}
\par\q Let $\Bbb N$ be the set of positive integers.
Let $a(x)=x+a_2x^2+a_3x^3+\cdots$ and
$\f{a(x)^m}{m!}=\sum_{n=m}^{\infty}a(n,m)\f{x^n}{n!}$ for $m\in\Bbb
N$. In Section 2 we show that for any $k,n\in\Bbb N$,
$$a(n+k,n)=\sum_{r=1}^k\bi{k-n}{k-r}\bi{k+n}{k+r}a(k+r,r).$$

\par Let $f(x)=c_0+c_1x+c_2x^2+\cdots$ with $c_0\not=0$.
 In Section 3 we
establish the following general inversion formula:
$$\align &a_n=n\sum_{m=1}^n[x^{n-m}]f(x)^m\cdot b_m\q(n=1,2,3,\ldots)
\\&\iff b_n=\f 1n\sum_{m=1}^n[x^{n-m}]f(x)^{-n}\cdot a_m\q(n=1,2,3,\ldots),
\endalign$$
 where $[x^k]g(x)$ is the coefficient of $x^k$ in
the power series expansion of $g(x)$. As a consequence, for a given
complex number $t$ we have the following inversion formula:
$$a_n=n\sum_{m=1}^n\bi{mt}{n-m}b_m\q (n\ge 1)
\iff b_n=\f 1n\sum_{m=1}^{n}\bi{-nt}{n-m}a_m\q (n\ge 1).$$

\par Let $\a^{-1}(x)$ be the inverse function of $\a(x)$. In Section 4
  we derive a general formula for
$[x^{m+n}]\a(x)^m$ by using the power series expansion of
$\a^{-1}(x)$. As a consequence, we deduce a symmetric inversion
formula, see Theorem 4.3.
\par Suppose $n\in\Bbb N$ and $k\in\{0,1,\ldots,n\}$.
Let $s(n,k)$ be the unsigned Stirling number of the first kind and
$S(n,k)$ be the Stirling number of the second kind defined by
$$x(x-1)\cdots(x-n+1)=\sum_{k=0}^n(-1)^{n-k}s(n,k)x^k$$ and
$$x^n=\sum_{k=0}^nS(n,k)x(x-1)\cdots(x-k+1).$$
In the paper we obtain new formulas for Stirling numbers, see
Theorems 2.3 and 4.2.

 \subheading{2. The formula for $[x^m]f(x)^t$}
 \pro{Lemma 2.1}
Let $t$ be a variable and $m\in\Bbb N$. Then
$$ \align &[x^m](1+a_1x+a_2x^2+\cdots+a_mx^m+\cdots)^t\\&=
\sum_{k_1+2k_2+\cdots+mk_m=m}\f {t(t-1)\cdots (t-(k_1+\cdots
+k_m)+1)}{k_1!\cdots k_m!}a_1^{k_1} \cdots a_m^{k_m}.\endalign $$
\endpro
Proof. Using the binomial theorem and the multinomial theorem we
see that
$$\aligned &[x^m](1+a_1x+a_2x^2+\cdots+a_mx^m+\cdots)^t\\&=
[x^m](1+a_1x+a_2x^2+\cdots+a_mx^m)^t\\&= [x^m]\sum_{n=0}^{\infty}\bi
t{n}(a_1x+a_2x^2+\cdots+a_mx^m)^n\\&=\sum_{n=0}^{m}\bi
t{n}[x^m](a_1x+a_2x^2+\cdots+a_mx^m)^n\\&=\sum_{n=0}^{m}\bi
tn[x^m]\sum_{k_1+k_2+\cdots+k_m=n}\f {n!}{k_1!\cdots
k_m!}(a_1x)^{k_1}\cdots (a_mx^m)^{k_m}\\&= \sum_{n=0}^m \bi
tn\sum\Sb k_1+\cdots+k_m=n\\k_1+2k_2+\cdots+mk_m=m\endSb \f
{n!}{k_1!\cdots k_m!}a_1^{k_1}\cdots
a_m^{k_m}\\&=\sum_{k_1+2k_2+\cdots+mk_m=m}\f {t(t-1)\cdots
(t-(k_1+\cdots +k_m)+1)}{k_1!\cdots k_m!}a_1^{k_1} \cdots a_m^{k_m}.
\endaligned $$ \pro{Theorem 2.1} Let $t$ be a variable, $m\in\Bbb N$
and $f(x)=1+a_1x+a_2x^2+\cdots$. Then
$$[x^m]f(x)^t=\sum_{r=1}^m\bi{m-t}{m-r}\bi tr[x^m]f(x)^r. $$
\endpro
Proof. From Lemma 2.1 we see that $[x^m]f(x)^t$ is a polynomial of
$t$ with degree $\le m$. Hence
$$P_m(t)=[x^m]f(x)^t-
\sum_{r=1}^m\bi{m-t}{m-r}\bi tr[x^m]f(x)^r $$ is also a polynomial
of $t$ with degree $\le m$. If $r\in\{1,2,\ldots,m\}$ and
$t\in\{0,1,\ldots,m\}$ with $t\not=r$, then $t<r$ or $m-t<m-r$ and
hence $\bi{m-t}{m-r}\bi tr=0$. Thus $P_m(t)=0$ for $t=0,1,\ldots,m$.
Therefore $P_m(t)=0$ for all $t$. This yields the result.
\pro{Corollary 2.1} Let $m\in\Bbb N$ and let $a$ be a complex
number. Then
$$\sum_{r=1}^m\bi{m+a}{m-r}(-1)^{m-r}\bi {a+r-1}rr^m=a^m.$$
\endpro
Proof. Clearly $[x^m](\t{e}^x)^t=\f{t^m}{m!}$. Thus, by Theorem 2.1
we have
$$\f{t^m}{m!}=\sum_{r=1}^m\bi{m-t}{m-r}\bi tr\f{r^m}{m!}.$$
Now taking $t=-a$ and noting that $\bi{-a}r=(-1)^r\bi{a+r-1}r$ we
deduce the result.

\pro{Theorem 2.2} Let $a(x)=x+a_2x^2+a_3x^3+\cdots$. For $m\in\Bbb
N$ let $\f{a(x)^m}{m!}=\sum_{n=m}^{\infty}a(n,m)\f{x^n}{n!}$. Then
for any $k,n\in\Bbb N$ we have
$$a(n+k,n)=\sum_{r=1}^k\bi{k-n}{k-r}\bi{k+n}{k+r}a(k+r,r).$$
\endpro
Proof. Set $\a(x)=a(x)/x$. Then for $m\in\Bbb N$ we have
$$\a(x)^m=\f{a(x)^m}{x^m}=\sum_{k=0}^{\infty}a(m+k,m)\cdot\f{m!}{(m+k)!}x^k.$$
Thus,
$$[x^k]\a(x)^n=a(n+k,n)\f{n!}{(n+k)!}
\qtq{and}[x^k]\a(x)^r=a(k+r,r)\f{r!}{(k+r)!}.$$ Since $\a(0)=1$, by
Theorem 2.1 we have
$$[x^k]\a(x)^n=\sum_{r=1}^k\bi{k-n}{k-r}\bi nr[x^k]\a(x)^r.$$
Hence
$$a(n+k,n)\f{n!}{(n+k)!}=\sum_{r=1}^k\bi{k-n}{k-r}\bi
nr\f{r!}{(k+r)!}a(k+r,r)$$ and so
$$a(n+k,n)=\sum_{r=1}^k\bi{k-n}{k-r}\f{(n+k)!}{(n-r)!(k+r)!}a(k+r,r).$$
This is the result.

\pro{Theorem 2.3} Let $k,n\in\Bbb N$. Then
$$S(n+k,n)=\sum_{r=1}^k\bi{k-n}{k-r}\bi{k+n}{k+r}S(k+r,r)$$
and
$$s(n+k,n)=\sum_{r=1}^k\bi{k-n}{k-r}\bi{k+n}{k+r}s(k+r,r).$$
\endpro
Proof. It is well known that ([1])
$$\f{(\t{e}^x-1)^m}{m!}=\sum_{n=m}^{\infty}S(n,m)\f{x^n}{n!}
\q \t{and}\q
\f{(\log(1+x))^m}{m!}=\sum_{n=m}^{\infty}(-1)^{n-m}s(n,m)\f{x^n}{n!}
.$$
  Thus the result follows from Theorem 2.2.
  \subheading{3. A
general inversion formula involving $[x^k]f(x)^t$}
 \pro{Lemma
3.1} Let $\a^{-1}(x)$ be the inverse function of $\a(x)$. Then for
any two sequences $\{a_n\}_{n=0}^{\infty}$ and
$\{b_n\}_{n=0}^{\infty}$ we have:
$$\align &a_n=\sum_{m=0}^{\infty}[x^n]\a(x)^mb_m\q(n=0,1,2,\ldots)
\\&\iff b_n=\sum_{m=0}^{\infty}[x^n]\a^{-1}(x)^ma_m\q(n=0,1,2,\ldots).
\endalign$$\endpro
Proof. Let $a(x)=\sum_{n=0}^{\infty}a_nx^n$ and
$b(x)=\sum_{n=0}^{\infty}b_nx^n$. Then clearly
$$\align
&a_n=\sum_{m=0}^{\infty}[x^n]\a(x)^mb_m\q(n=0,1,2,\ldots)
\\&\iff
a(x)=\sum_{m=0}^{\infty}b_m\sum_{n=0}^{\infty}[x^n]\a(x)^mx^n
=\sum_{m=0}^{\infty}b_m\a(x)^m
\\&\iff a(x)=b(\a(x))
\iff b(x)=a(\a^{-1}(x))
\\&\iff b_n=\sum_{m=0}^{\infty}[x^n]\a^{-1}(x)^ma_m\q(n=0,1,2,\ldots).
\endalign$$ So the lemma is proved.
\par\q
\pro{Theorem 3.1} Let $k\in\Bbb N$. For nonnegative integers $m$ and
$n$ let
$$\a_k(n,m)=\cases (-1)^{\f nk}\bi{\f mk}{\f nk}&\t{if $k\mid n$,}
\\0&\t{if $k\nmid n$.}\endcases$$
Then we have the following inversion formula:
$$\aligned &a_n=\sum_{m=0}^{\infty}\a_k(n,m)b_m\q(n=0,1,2,\ldots)
\\&\iff b_n=\sum_{m=0}^{\infty}\a_k(n,m)a_m\q(n=0,1,2,\ldots).
\endaligned$$\endpro
Proof. Let $\a(x)=(1-x^k)^{\f 1k}\ (0<x<1)$. Then clearly
$\a^{-1}(x)=\a(x)$ and $\a(x)^m=(1-x^k)^{\f mk}=\sum_{r=0}^{\infty}
\bi{\f mk}r(-1)^rx^{kr}=\sum_{n=0}^{\infty}\a_k(n,m)x^n$. Thus
applying Lemma 3.1 we deduce the theorem.
 \pro{Lemma 3.2 (Lagrange
inversion formula ([1, p.148], [3, pp.36-44]))} \newline Let
$\a(x)=\a_1x+\a_2x^2+\cdots$ with $\a_1\not=0$, and let $k,n\in\Bbb
N$ with $k\le n$. Then
$$[x^n](\a^{-1}(x))^k=\f kn[x^{n-k}]\Ls{\a(x)}x^{-n}.$$\endpro
\pro{Theorem 3.2} Let $f(x)=c_0+c_1x+c_2x^2+\cdots$ with
$c_0\not=0$. Then for any two sequences $\{a_n\}$ and $\{b_n\}$ we
have the following inversion formula:
$$\align &a_n=n\sum_{m=1}^n[x^{n-m}]f(x)^m\cdot b_m\q(n=1,2,3,\ldots)
\\&\iff b_n=\f 1n\sum_{m=1}^n[x^{n-m}]f(x)^{-n}\cdot a_m\q(n=1,2,3,\ldots).
\endalign$$\endpro
Proof. Set $\a(x)=xf(x)$. Then clearly $[x^n]\a(x)^m=0$ for $m>n$.
As $\a^{-1}(xf(x))=\a^{-1}(\a(x))=x$ we see that $\a^{-1}(0)=0$ and
so $\a^{-1}(x)=d_1x+d_2x^2+\cdots$ for some $d_1,d_2,\ldots$. Thus
$[x^n]\a^{-1}(x)^m=0$ for $m>n$. Set $a_0=b_0=0$. From Lemma 3.1 we
see that
$$\align &a_n=\sum_{m=0}^{\infty}[x^n]\a(x)^m\cdot b_m
=\sum_{m=1}^n[x^n]\a(x)^m\cdot b_m\q(n\ge 1)
\\&\iff b_n=\sum_{m=0}^{\infty}[x^n]\a^{-1}(x)^m\cdot a_m
=\sum_{m=1}^n[x^n]\a^{-1}(x)^m\cdot a_m\q(n\ge 1).
\endalign$$ For $m\le n$ we see that
$[x^n]\a(x)^m=[x^n]x^mf(x)^m=[x^{n-m}]f(x)^m$ and
$[x^n]\a^{-1}(x)^m=\f mn[x^{n-m}]f(x)^{-n}$ by Lemma 3.2. Thus
$$a_n=\sum_{m=1}^n[x^{n-m}]f(x)^m\cdot b_m\ (n\ge 1) \iff
b_n=\sum_{m=1}^n\f mn[x^{n-m}]f(x)^{-n}\cdot a_m\ (n\ge 1).$$ Now
substituting $a_n$ by $a_n/n$ we obtain the result.
\par\q
\par As $\t{e}^{cx}=\sum_{k=0}^{\infty}\f{(cx)^k}{k!}$, we see that
$$[x^{n-m}](\t{e}^x)^m=\f{m^{n-m}}{(n-m)!}\qtq{and}
[x^{n-m}](\t{e}^x)^{-n}=\f{(-n)^{n-m}}{(n-m)!}.$$ Thus, putting
$f(x)=\t{e}^x$ in Theorem 3.2 we have the following inversion
formula:
$$a_n=n\sum_{m=1}^n\f{m^{n-m}}{(n-m)!}b_m\q(n\ge 1)
\iff b_n=\f 1n\sum_{m=1}^n\f{(-n)^{n-m}}{(n-m)!}a_m\q(n\ge 1).$$
Substituting $a_n$ by $a_n/(n-1)!$, and $b_n$ by $b_n/n!$ we obtain
$$a_n=\sum_{m=1}^n\bi nm m^{n-m}b_m\q(n\ge 1)
\iff b_n=\sum_{m=1}^n\bi{n-1}{m-1}(-n)^{n-m}a_m\q(n\ge 1).$$ This is
a known result. See [2, p.96].
\par As $(1+x)^{ct}=\sum_{k=0}^{\infty}\bi{ct}kx^k\ (|x|<1)$, we
see that
$$[x^{n-m}](1+x)^{mt}=\bi{mt}{n-m}\qtq{and}[x^{n-m}](1+x)^{-nt}
=\bi{-nt}{n-m}.$$ Now putting $f(x)=(1+x)^t$ in Theorem 3.2 and
applying the above we deduce the following result. \pro{Theorem 3.3}
Let t be a complex number. For any two sequences $\{a_n\}$ and
$\{b_n\}$ we have the following inversion formula:
$$a_n=n\sum_{m=1}^n\bi{mt}{n-m}b_m\q (n\ge 1)
\iff b_n=\f 1n\sum_{m=1}^{n}\bi{-nt}{n-m}a_m\q (n\ge 1).$$
 \endpro

 \pro{Theorem 3.4} Let
$f(x)=c_0+c_1x+c_2x^2+\cdots$ with $c_0\not=0$. For $k,n\in\Bbb N$
with $k<n$ we have
$$\sum_{m=k}^n\f
1m[x^{n-m}]f(x)^m\cdot [x^{m-k}]f(x)^{-m}=\sum_{m=k}^n
m[x^{m-k}]f(x)^k\cdot [x^{n-m}]f(x)^{-n}=0.$$
\endpro
Proof. For $m\in\Bbb N$ let $b_m=\f
1m\sum_{k=1}^m[x^{m-k}]f(x)^{-m}\cdot y^k$. Applying Theorem 3.2 we
see that
$$ \sum_{m=1}^n[x^{n-m}]f(x)^m\cdot b_m=\f {y^n}n.$$
On the other hand,
$$\align \sum_{m=1}^n[x^{n-m}]f(x)^m\cdot b_m
&=\sum_{m=1}^n[x^{n-m}]f(x)^m\cdot \f
1m\sum_{k=1}^m[x^{m-k}]f(x)^{-m}\cdot y^k
\\&=\sum_{k=1}^n\Big(\sum_{m=k}^n[x^{n-m}]f(x)^m\cdot \f 1m[x^{m-k}]f(x)^{-m}\Big)y^k.
\endalign$$ Thus,
$$\sum_{k=1}^n\Big(\sum_{m=k}^n\f 1m[x^{n-m}]f(x)^m[x^{m-k}]f(x)^{-m}\Big)y^k=\f{y^n}n$$
and hence
$$\sum_{m=k}^n\f
1m[x^{n-m}]f(x)^m\cdot [x^{m-k}]f(x)^{-m}=0\qtq{for}k<n.$$
\par
For $m\in\Bbb N$ let $a_m=m\sum_{k=1}^m[x^{m-k}]f(x)^k\cdot y^k$.
Applying Theorem 3.2 we have
$$ \sum_{m=1}^n[x^{n-m}]f(x)^{-n}\cdot a_m=ny^n.$$
On the other hand,
$$\align \sum_{m=1}^n[x^{n-m}]f(x)^{-n}\cdot a_m
&=\sum_{m=1}^n[x^{n-m}]f(x)^{-n}\cdot
m\sum_{k=1}^m[x^{m-k}]f(x)^k\cdot y^k
\\&=\sum_{k=1}^n\Big(\sum_{m=k}^n[x^{n-m}]f(x)^{-n}\cdot m[x^{m-k}]f(x)^k\Big)y^k.
\endalign$$ Thus,
$$\sum_{k=1}^n\Big(\sum_{m=k}^nm[x^{m-k}]f(x)^k\cdot [x^{n-m}]f(x)^{-n} \Big)y^k=ny^n$$
and hence
$$\sum_{m=k}^nm[x^{m-k}]f(x)^k\cdot [x^{n-m}]f(x)^{-n}=0\qtq{for}k<n.$$
This completes the proof.
 \pro{Corollary 3.1} For $k,n\in\Bbb N$ with $k<n$ we have
$$\sum_{m=k}^n\f
1m\binom{mt}{n-m}\binom{-mt}{m-k}=\sum_{m=k}^n
m\binom{kt}{m-k}\binom{-nt}{n-m}=0.$$
\endpro
Proof. Since $(1+x)^{rt}=\sum_{s=0}^{\infty}\binom{rt}sx^s$, taking
$f(x)=(1+x)^t$ in Theorem 3.4 we deduce the result.
 \subheading{4. A formula for $[x^{m+n}]\a(x)^m$}
 \pro{Theorem 4.1}
Let $\b(x)=x\sum_{n=0}^{\infty}\b_nx^n$ with $\b_0\neq0$. Let
$\a(x)$ be the inverse function of $\b(x)$. For $m,n\in\Bbb N$ we
have
$$\align [x^{m+n}]\a(x)^m&=\f m{(n+m)!}
\sum_{\sn=n}\f{(n+m-1+k_1+\cdots+k_n)!}{\kn1}\\&\q\times
(-1)^{k_1+k_2+\cdots+k_n}\b_0^{-n-m-k_1-\cdots-k_n}\b_1^{k_1}\b_2^{k_2}
\cdots\b_n^{k_n}.\endalign$$\endpro Proof. By the multinomial
theorem we have
$$\Big(\sum_{k=1}^n\f{\b_k}{\b_0}x^k\Big)^s=\sum_{k_1+\cdots+k_n=s}\f{s!}
{k_1!\cdots k_n!}\prod_{i=1}^n\Big(\f{\b_i}{\b_0}x^i\Big)^{k_i}.$$
Thus
$$[x^n]\Big(\sum_{k=1}^{\infty}
\f{\b_k}{\b_0}x^k\Big)^s=[x^n]\Big(\sum_{k=1}^n\f{\b_k}{\b_0}x^k\Big)^s
=\sum\Sb k_1+\cdots+k_n=s\\k_1+2k_2+\cdots+nk_n=n\endSb\f{s!}
{k_1!\cdots k_n!}\prod_{i=1}^n\Big(\f{\b_i}{\b_0}\Big)^{k_i}.$$ As
$$\aligned &\b_0^{m+n}\Ls x{\b(x)}^{m+n}-1\\&
=\b_0^{m+n}\Big(\b_0+\sum_{k=1}^{\infty}\b_kx^k\Big)^{-n-m}-1
=\Big(1+\sum_{k=1}^{\infty}\f{\b_k}{\b_0}x^k\Big)^{-n-m}-1
\\&=\sum_{s=1}^{\infty}\f{(-n-m)(-n-m-1)\cdots(-n-m-s+1)}{s!}
\Big(\sum_{k=1}^{\infty}\f{\b_k}{\b_0}x^k\Big)^s.
\endaligned$$ From the above
we see that
$$\aligned &[x^n]\b_0^{m+n}\Ls x{\b(x)}^{m+n}\\&
=\sum_{s=1}^{\infty}\f{(-n-m)(-n-m-1)\cdots(-n-m-s+1)}{s!}
\\&\qq\times\sum\Sb k_1+\cdots+k_n=s\\k_1+2k_2+\cdots+nk_n=n\endSb\f{s!}
{k_1!\cdots k_n!}\prod_{i=1}^n\Big(\f{\b_i}{\b_0}\Big)^{k_i}
\\&=\sum_{k_1+2k_2+\cdots+nk_n=n}\f{(n+m)(n+m+1)\cdots
(n+m+k_1+\cdots+k_n-1)}{k_1!\cdots k_n!}\\&\qq\times\Big(-\f
1{\b_0}\Big)^{k_1+\cdots+k_n}\b_1^{k_1}\cdots\b_n^{k_n}.
\endaligned$$
Thus applying Lemma 3.2 we have
$$\aligned [x^{m+n}]\a(x)^m&=\f m{n+m}[x^n]\Ls x{\b(x)}^{m+n}
\\&=\f m{n+m}\b_0^{-m-n}\sum_{k_1+2k_2+\cdots+nk_n=n}\f{(k_1+\cdots+k_n+n+m-1)!}{
k_1!\cdots
k_n!(n+m-1)!}\\&\qq\times(-1)^{k_1+\cdots+k_n}\b_0^{-(k_1+\cdots+k_n)}
\b_1^{k_1}\cdots\b_n^{k_n}.
\endaligned.$$ This yields the result.
\pro{Corollary 4.1} For $m,n\in\Bbb N$ we have
$$\sum_{\sn=n}\f{(k_1+\cdots+k_n+m+n-1)!}{(m+n-1)!\kn1}(-1)^{k_1+\cdots+k_n}
=(-1)^n\binom{m+n}m. $$\endpro Proof. Let
$\b(x)=x\sum_{r=0}^{\infty}x^r=\f x{1-x}$. Then the inverse function
of $\b(x)$ is given by $\a(x)=\f x{1+x}$. Using the binomial theorem
we see that $[x^{m+n}]\a(x)^m=[x^n](1+x)^{-m}
=\binom{-m}n=(-1)^n\binom{m+n-1}n$. Now applying Theorem 4.1 we
deduce the result.

\pro{Corollary 4.2} For $n\in\Bbb N$ we have
$$\aligned&\sum_{\sn=n}\f{(k_1+\cdots+k_n+n)!}{\kn1}(-1)^{k_1+\cdots+k_n}
 2^{k_1}3^{k_2}\cdots (n+1)^{k_n}\\&=(-1)^n\cdot(n+1)!\cdot\f
1{n+2}\bi {2n+2}{n+1}.\endaligned$$\endpro
 Proof. Let
$$\b(x)=\f x{(1+x)^2}\qtq{and}\a(x)=\f{1-\sqrt{1-4x}}{2x}-1\
 \big(0<x<\f 14\big).$$ It is easily seen that $\a(x)=\b^{-1}(x)$. From
 the binomial theorem we know that
$$\a(x)=x
\sum _{n=0}^{\infty}\f 1{n+2}\bi {2n+2}{n+1}x^n \qtq{and}\b(x)
=x\sum_{n=0}^{\infty}(-1)^n(n+1)x^n.$$ Now applying Theorem 4.1
(with $m=1$) we deduce the result.
\pro{Corollary 4.3} For $n\in\Bbb
N$ we have
$$\sum_{\sn=n}\f{(k_1+\cdots+k_n+2n)!}{\kn1}\cdot
\f{(-1)^{k_1+k_2+\cdots+k_n+n}}
 {3!^{k_1}5!^{k_2}\cdots (2n+1)!^{k_n}}=(2n-1)!!^2.$$\endpro
Proof. It is well known that
$$\sin
x=\sum_{n=0}^{\infty}\f{(-1)^n}{(2n+1)!}x^{2n+1}$$ and $$\arcsin
x=x+\sum_{n=1}^{\infty}\f{(2n-1)!!}{(2n+1)\cdot(2n)!!}x^{2n+1}\
(|x\le 1).$$ Set $\b(x)=\sin x=x\sum_{n=0}^{\infty}\b_nx^n$. Then
$\b^{-1}(x)=\arcsin x$ and
$$\b_i=\cases 0&\t{if $2\nmid i$,}
\\\f{(-1)^{i/2}}{(i+1)!}&\t{if $2\mid i$.}\endcases$$
Thus, taking $m=1$ in Theorem 4.1 and substituting $n$ by $2n$ we
obtain
$$\align &(2n+1)!\cdot [x^{2n+1}]\arcsin x\\&=
\sum_{k_1+2k_2+\cdots+2nk_{2n}=2n}\f{(2n+k_1+k_2\cdots+k_{2n})!}{k_1!k_2\cdots
k_{2n}!}(-1)^{k_1+k_2\cdots+k_{2n}}\b_1^{k_1}\b_2^{k_2}\cdots
\b_{2n}^{k_{2n}}
\\&=\sum_{k_2+2k_4+\cdots+nk_{2n}=n}\f{(2n+k_2+k_4\cdots+k_{2n})!}{k_2!k_4!\cdots
k_{2n}!}(-1)^{k_2+k_4+\cdots+k_{2n}}\prod_{i=1}^n\Ls{(-1)^i}{(2i+1)!}^{k_{2i}}.
\endalign$$
Replacing $k_{2i}$ by $k_i$ in the above formula and observing that
$$(2n+1)!\cdot [x^{2n+1}]\arcsin x=(2n+1)!\cdot
\f{(2n-1)!!}{(2n+1)\cdot (2n)!!}=(2n-1)!!^2$$ we deduce the result.

 \pro{Theorem 4.2} For $m,n\in\Bbb N$ we have
$$S(n+m,m)=\f {(-1)^n}{(m-1)!}\sum_{k_1+2k_2+\cdots+nk_n=n}(-1)^{k_1+\cdots+k_n}
 \f{(k_1+\cdots+k_n+n+m-1)!}{2^{k_1} k_1!\cdot 3^{k_2}
k_2!\cdots (n+1)^{k_n} k_n!}$$ and
$$s(n+m,m)=\f{(-1)^n}{(m-1)!}\sum_{k_1+2k_2+\cdots+nk_n=n}(-1)^{k_1+\cdots+k_n}
 \f{(k_1+\cdots+k_n+n+m-1)!}{2!^{k_1} k_1!\cdot 3!^{k_2}
k_2!\cdots (n+1)!^{k_n}k_n!}.$$
\endpro
Proof. Clearly $\t{e}^x-1$ and $\log(1+x)$ are a pair of inverse
functions. As
$$\f{(\t{e}^x-1)^m}{m!}=\sum_{n=0}^{\infty}S(n+m,m)\f{x^{n+m}}{(n+m)!}
\qtq{and}\log(1+x)=\sum_{i=0}^{\infty}\f{(-1)^i}{i+1}x^{i+1},$$
putting $\a(x)=\t{e}^x-1$, $\b(x)=\log(1+x)$ and
 $\b_i=\f{(-1)^i}{i+1}$ in Theorem 4.1
we see that
$$\align \f{m!S(n+m,m)}{(n+m)!}&
=[x^{m+n}](\t{e}^x-1)^m
\\&=\f m{(n+m)!}\sum_{k_1+2k_2+\cdots+nk_n=n}
 \f{(k_1+\cdots+k_n+n+m-1)!}{k_1!\cdots k_n!}
 \\&\qq\times(-1)^{k_1+k_2+\cdots+k_n}\cdot(-1)^{k_1+2k_2+\cdots+nk_n}\f 1{2^{k_1}\cdot 3^{k_2}\cdots
 (n+1)^{k_n}}.\endalign$$
 Since
$$\f{(\log(1+x))^m}{m!}=\sum_{n=0}^{\infty}(-1)^{n}s(n+m,m)\f{x^{n+m}}{(n+m)!}
\qtq{and}\t{e}^x-1=\sum_{i=0}^{\infty}\f{x^{i+1}}{(i+1)!},$$
putting $\a(x)=\log(1+x)$, $\b(x)=\t{e}^x-1$ and
 $\b_i=\f 1{(i+1)!}$ in Theorem 4.1
we see that
$$\align &(-1)^n\f{m!s(n+m,m)}{(n+m)!}\\&
=[x^{m+n}](\log(1+x))^m
\\&=\f m{(n+m)!}\sum_{k_1+2k_2+\cdots+nk_n=n}
 \f{(k_1+\cdots+k_n+n+m-1)!}{k_1!\cdots k_n!}\cdot
 \f{(-1)^{k_1+\cdots+k_n}}{2!^{k_1}\cdot 3!^{k_2}\cdots
 (n+1)!^{k_n}}.\endalign$$
By the above, the theorem is proved.
\par
\par We remark that Theorem 4.2 provides a straightforward method to
calculate $s(n+m,m)$ and $S(n+m,m)$ for small $n$. For example, we
have
$$S(m+3,m)=\binom{m+1}2\binom{m+3}4\qtq{and}s(m+3,m)
=\binom{m+3}2\binom{m+3}4.\tag 4.1$$
\pro{Corollary 4.4} For
$m,n\in\Bbb N$ we have
$$\align &\sum_{r=0}^m\binom mr(-1)^{m-r}r^{m+n}
\\&=m\sum_{k_1+2k_2+\cdots+nk_n=n}(-1)^{k_1+\cdots+k_n+n}
 \f{(k_1+\cdots+k_n+n+m-1)!}{2^{k_1} k_1!\cdot 3^{k_2}
k_2!\cdots (n+1)^{k_n} k_n!}.\endalign$$
\endpro
Proof. It is well known that ([1, p.204])
$$\sum_{r=0}^m\binom
mr(-1)^{m-r}r^{m+n}=m!S(n+m,m).$$ Combining this with Theorem 4.2 we
obtain the result.
\par Let
$\a(x)=-x+\a_1x^2+\a_2x^3+\cdots$ and
$\b(x)=-x+\b_1x^2+\b_2x^3+\cdots$ be a pair of inverse functions.
Taking $m=1$ in Theorem 4.1 we deduce:
 \pro{Theorem 4.3} We have the following inversion formula:
$$ \align&\a_n=\f {(-1)^{n+1}}{(n+1)!}\sum _{\sn =n}\f {(k_1+\cdots
+k_n+n)!}{\kn1}\b_1^{k_1}\cdots \b_n^{k_n}(n\ge 1)
\\&\iff\b_n=\f {(-1)^{n+1}} {(n+1)!}\sum_{\sn=n} \f {(k_1+\cdots
+k_n+n)!} {\kn1} \a_1^{k_1}\cdots \a_n^{k_n}(n\ge 1).
\endalign$$\endpro
 \pro{Definition
4.1} If $\a(x)=\a^{-1}(x)$, we say that $\a(x)$ is a self-inverse
function.\endpro
\par For example, $\a(x)=\f{rx+s}{tx-r}\ ((r^2+t^2)(r^2+st)\not=0)$ and
$\a(x)=(1-x^k)^{\f 1k}$ are self-inverse functions.

 \pro{Theorem 4.4} Let
$\a(x)=-x+\a_1x^2+\a_2x^3+\cdots$ be a self-inverse function. Then
$\a_2,\a_4,\ldots$ depend only on $\a_1,\a_3,\ldots$. Moreover, for
$n\in\Bbb N$,
$$\aligned &\sum_{k_1+2k_2+\cdots+(n-1)k_{n-1}=n}
\f {(k_1+\cdots+k_{n-1}+n)!}{k_1!\cdots k_{n-1}!} \a_1^{k_1}\cdots
\a_{n-1}^{k_{n-1}}\\&\qq\qq=\cases 0&\t{if $2\nmid n$,}
\\-2\cdot (n+1)!\a_n&\t{if $2\mid n$.}\endcases\endaligned\tag 4.2$$
\endpro
Proof. By Theorem 4.3 we have
$$\align\a_n&=\f
{(-1)^{n+1}}{(n+1)!}\sum_{\sn=n}\f {(k_1+\cdots+k_n+n)!}{\kn1
}\a_1^{k_1}\cdots\a_{n}^{k_{n}}\\&=\f
{(-1)^{n+1}}{(n+1)!}\sum_{k_1+2k_2+\cdots+(n-1)k_{n-1}=n}\f
{(k_1+\cdots+k_{n-1}+n)!}{k_1!\cdots k_{n-1}!
}\a_1^{k_1}\cdots\a_{n-1}^{k_{n-1}}+(-1)^{n+1}\a_n.\endalign$$ Thus
(4.2) is true.  Using (4.2) and induction we deduce that
$\a_2,\a_4,\ldots$ depend only on $\a_1,\a_3,\ldots$. This completes
the proof.
\par If $\a(x)=-x+\a_1x^2+\a_2x^3+\cdots$ is a self-inverse
function, from (4.2) we deduce
$$\aligned&\a_2=-\a_1^2,\ \a_4=2\a_1^4-3\a_1\a_3,
\\&\a_6=-13\a_1^6-4\a_1\a_5-2\a_3^2+18\a_1^3\a_3,
\\&\a_8=145\a_1^8-221\a_1^5\a_3+50\a_1^2\a_3^2+35\a_1^3\a_5-5\a_3\a_5-5\a_1\a_7.
\endaligned\tag 4.3$$

\enddocument
\bye